\documentclass{article}
\usepackage{tikz}
\usetikzlibrary{arrows,shapes,automata,petri,positioning,calc}
\usepackage{amsmath}
\usepackage{amsthm}
\theoremstyle{definition}
\newtheorem{definition}{Definition}

\usepackage{graphicx}
\usepackage{wrapfig}
\usepackage{caption}
\usepackage{csquotes}
\usepackage{subcaption}
\usepackage{mathtools}
\usepackage{sidecap}
\usepackage{hyperref}
\usepackage[title]{appendix}%
\usepackage{algorithm}
\usepackage{bbm}
\usepackage{amsfonts}
\usepackage{fancyhdr}
\usepackage{mathrsfs}
\usepackage{float}
\usepackage{enumerate}
\usepackage{enumitem}
\usepackage[english]{babel}
\usepackage[nottoc]{tocbibind}
\usepackage{authblk}

\usepackage[backend=biber,giveninits=true,bibstyle=authoryear,url = false,doi = false,isbn=false]{biblatex}
\DeclareNameAlias{sortname}{family-given}
\DeclareNameAlias{default}{family-given}

\DeclareFieldFormat[inbook]{title}{#1}
\DeclareFieldFormat[article]{title}{#1}
\DeclareFieldFormat[inproceedings]{title}{#1}

\defbibenvironment{bibliography}
  {\enumerate
     {}
     {\setlength{\labelwidth}{\labelnumberwidth}%
      \setlength{\leftmargin}{\labelwidth}%
      \setlength{\labelsep}{\biblabelsep}%
      \addtolength{\leftmargin}{\labelsep}%
      \setlength{\itemsep}{\bibitemsep}%
      \setlength{\parsep}{\bibparsep}}%
      }
  {\endlist}
  {\item}
\title{Minimising Numbers of Losses and Abandonments in Small Call Centres Under a Transient Regime}
\author[1,2]{Mark Fackrell\thanks{\href{mailto:fackrell@unimelb.edu.au}{fackrell@unimelb.edu.au}}}
\author[1,2]{Hritika Gupta\thanks{\href{mailto:hritikag@student.unimelb.edu.au}{hritikag@student.unimelb.edu.au}}}
\author[1,2]{Peter G. Taylor\thanks{\href{mailto:taylorpg@unimelb.edu.au}{taylorpg@unimelb.edu.au}}}
\affil[1]{School of Mathematics and Statistics, The University of Melbourne}
\affil[2]{ARC Training Centre in Optimisation Technologies, Integrated Methodologies, and Applications (OPTIMA)}

\begin{document}

\maketitle

\begin{abstract}
    We consider the problem of allocating customers to agents in small call centres so that transient performance indicators in terms of expected numbers of losses and abandonments are optimised. To gain insight into the general structure of optimal policies, we start by (i) using backward induction based upon Bellman's equation for a finite-horizon discrete-time model, and (ii) deriving stationary policies for an infinite-horizon continuous-time model with discounting. 
    
    We then compare the performance of such policies in our original finite-horizon continuous-time model using a first-step analysis applied to the Laplace transforms of the relevant measures. 
\end{abstract}

\section{Introduction}\label{intro}

The term \textit{call centre} is commonly used to describe a telephone-based human-service operation. A call centre provides teleservices in which customers and service agents are remote from each other (\textcite{mandelbaum2000empirical}). They have become increasingly popular, and service quality is now a crucial factor in attracting and retaining customers.

Call centre service providers frequently have contracts with their clients under which they have to manage their client's customer care department and meet certain quality standards, such as serving at least $90\%$ of the calls within 20 seconds or having less than $3\%$ of call abandonments per month. Human resources typically account for more than $60\%$ of the total operating cost of a call centre (\textcite{gans_telephone_2003}, \textcite{brown_statistical_2005}). Hence, a reasonable objective for a call centre service provider is to minimise the staffing requirements while meeting all the targets. This can be achieved by improving the prediction accuracy related to call arrival processes, using an efficient agent scheduling method and making the best use of resources by appropriately allocating incoming calls to agents. It is the third of these that we will focus on in this paper.


In many call centres, customer queries are divided into different categories depending on the level of complexity of the query. Similarly, agents are also assigned a skill level depending on their ability to solve these queries. When a call arrives, based on the caller's initial input to an Interacting Voice Response (IVR) system, they are assigned to an agent with the appropriate skill level. Some agents are multi-skilled and can deal with more than one type of query. When a call centre has multi-skilled agents, there can be a number of possibilities for allocating an incoming call to an agent. We call the rules of assigning an incoming call to an agent a \textit{policy}.

Skills-based routing has been studied many times using different techniques. 
\textcite{wallace_staffing_2005} developed a suboptimal routing algorithm for skills-based routing, and showed using simulation analysis that limited cross training among agents could be almost as efficient as training all the agents for all the skills.
\textcite{stolletz2004performance} calculated the steady state probabilities and performance measures for a given priority based routing policy. 
\textcite{bhulai2009dynamic} obtained nearly optimal routing policies by minimising the long term expected average cost for the specific case where there is either no buffer for waiting customers or there is a buffer and all agents have either a single skill or are fully cross-trained.
\textcite{mehrbod2018caller} and \textcite{jorge2020intelligent} used machine learning techniques such as Random Forests and Neural Networks for their analysis. Simulation, see for example \textcite{adetunji2007performance}, \textcite{akhtar2010exploiting} and \textcite{mehrotra2012routing}, has also been a popular method for comparing different routing policies.


Our study is motivated by a problem brought to us by an industry partner that manages call centres for its clients. The company was struggling to meet targets in terms of expected numbers of abandonments and losses, despite having high non-occupancy rates among its agents. It wanted us to look at a particular client's call centre that was reasonably small, with around a dozen agents in total. There were four levels of agents of increasing expertise who were required to handle four levels of queries of increasing complexity. The company had good data on the numbers of arrivals of the various levels of query in half-hour intervals. This data demonstrated that, while the arrival processes could reasonably be modelled as homogeneous over a half-hourly time-scale, they were inhomogeneous over longer time scales.

While our industry partner utilised models to estimate its workload and determine the required total number of agents, it did not use mathematically-backed strategies for allocating calls to agents. Indeed, call allocation was managed by an employee who manually decided how to allocate calls based on observing the system.  We thought that a natural analysis would be to use a sequence of finite-horizon continuous-time Markov decision models for each half-hour interval. Since it is unreasonable that the state at the end of one interval would follow the stationary distribution of the next, such an analysis would necessarily have to be transient. 


However, transient analysis of finite-horizon continuous-time Markov decision processes is difficult to carry out because the number of transitions occurring within the time horizon is both random and policy-dependent (\textcite{van2018uniformization}). 
Furthermore, in many models, including ours, decisions need to be made only when particular events occur. In our context, these events are call arrivals or service completions in which a choice needs to be made about which agent, if any, to allocate to an arriving or waiting customer. 
 

A widely used technique for the transient analysis of continuous-time Markov chains is uniformisation. For a continuous-time Markov chain with transition matrix $Q$ that has bounded diagonal entries, this technique essentially writes the matrix exponential in the form
\begin{equation}
    \label{eq:uniformisation}
e^{Qt} = e^{-\lambda t}\sum_{k=0}^\infty \frac{(\lambda t)^k}{k!} P^k
\end{equation}
where $\lambda \geq \sup_{i} |q_{ii}|$ and $P=Q/\lambda + I$ with $I$ the identity matrix (see \textcite{gross1984randomization}, \textcite{buchholz2011numerical}).

The matrix $P$ on the right hand side of {\eqref{eq:uniformisation} is the transition matrix of a discrete-time Markov chain, which opens the door to a discrete-time analysis.
Uniformisation can be incorporated into a Markov decision process by letting the entries of $P$ depend on state-dependent actions. While this approach can be effective for many problems, it is not the most efficient approach for our problem. Observe the Markov chain with transition matrix $P$ frequently has states with a high probability of `self-transition' back to the same state. Putting that together with the above observation that, in our setting, decisions are required only at certain transition points, applying uniformisation would introduce a large number of redundant decision points, significantly increasing the computational burden. Since our goal is to develop an optimal policy that is robust, scalable, and computationally practical for a real-world call centre with multi-level agent hierarchies, we decided not to employ uniformisation.

The approach that we discuss below is related to the framework presented in \textcite[Sections 10.5 and 10.6]{kulkarni2014modeling}. Kulkarni discussed semi-Markov decision processes (SMDPs), where decisions are made at transition epochs occurring at random times, and the system evolves according to continuous-time dynamics between transitions. Our model is Markovian, but this approach still works.

We first analyse
\begin{enumerate}
\item the discrete-time jump chain embedded in the continuous time process over a fixed finite number of transitions and
\item an approximating continuous-time Markov decision process over an infinite horizon with discounting.
\end{enumerate}
For the first model, we make use of backward induction using Bellman's equation to characterise the value function and derive an optimal policy. In doing this, we have ignored the variability in the sojourn times that occurs in the original continuous-time model. A key question that remains is how to match the bound on the number of transitions in the discrete-time model to a finite-horizon analysis in continuous time. 

Our analysis of the second model is necessarily stationary. However, again we can expect to gain insight into the structure of policies that are optimal for the original model. 


Based on the hierarchical structure of our industry partner's centre, we first consider a model with $\nu$ classes of query and $\nu$ levels of agents, where level $i$ agents have the skills to deal with level $i$ and level $i-1$ queries. We work with this model in Sections \ref{General}, \ref{ch:res}, and \ref{ab4}.

To illustrate the approach, in Section \ref{General} we take $\nu = 2$ and use Bellman's equation to find the optimal policy for the embedded discrete-time Markov decision process over a finite horizon and the discounted continuous-time Markov decision process over an infinite horizon. This gives us an indication of the class of policies that are optimal. We discuss this in Section \ref{ch:res}. In Section \ref{ab4}, we revert back to the model of our industry partner's centre that has $\nu = 4$, adapting the method of \textcite{chiera_what_2002} to calculate the exact expected number of abandonments for the continuous-time Markov chain for a class of policies motivated by our earlier analysis. 

We then present a discussion of the computational aspects and challenges of this method in Section \ref{dis}. The paper concludes with Section \ref{con} where we summarise our findings, and present some ideas for future research.

\section{The optimal policy for $\nu=2$. }\label{General}

In this section, we work with a downsized problem. Consider a call centre with two levels of callers and two levels of agents, where level 1 agents can take only level 1 calls and level 2 agents can serve both level 1 and 2 calls. We aim to compute the policy that minimises the total expected cost of abandonments for a given number of agents. This optimal policy is fully flexible in the sense that it takes the optimal state-dependent action whenever there is an arrival or service without assuming that there is any predefined functional form of the policy.

\subsection{Events requiring decision making}\label{DM}

In the call centre system, two types of events result in an opportunity to allocate a customer to an unoccupied server: when there is a call arrival, and when there is a service completion. When a call arrives, we need to decide which available agent to assign that call to, or we can choose to add the call to a queue. When an agent serves a call, we need to decide which call out of those waiting (if any) we should assign to this agent, or we can choose to reserve the agent for another call that might arrive later.

In order to reduce our action space, we make the assumption that if a level 1 call arrives or is in the queue and a level 1 agent is available, the optimal action is to assign the level 1 call to the available level 1 agent. This makes sense, since level 1 agents are dedicated to level 1 calls and there is nothing to be gained by making a level 1 customer wait when such an agent is available. 

When an arrival occurs, we need to make a decision only when
\begin{itemize}
\item the arrival is a level 1 call, all level 1 agents are busy and level 2 agents are available, or 
\item the arrival is a level 2 call and level 2 agents are available.
\end{itemize}
The possible actions in these cases are to assign the arriving call to the available level 2 agent or not. We are not ruling out the case that level 2 agents might be reserved for level 1 calls, and hence might not take an arriving level 2 call even when they are available, although we expect that this will not happen for most sensible parameter values.

Similarly, when there is a service, we need to make a decision only when 
\begin{itemize}
\item the service is completed by a level 2 agent and either only level 1 or only level 2 calls are waiting in the queue, or
\item the service is completed by a level 2 agent and both level 1 and level 2 calls are waiting in the queue.
\end{itemize}
The possible actions in the first case are to assign the waiting call to the level 2 agent who just completed the service or not, and the possible actions in the second case are to assign a level 1 or a level 2 call to the level 2 agent who just completed the service. A level 2 agent cannot reject both level 1 and level 2 calls. 

\subsection{The model}\label{gen_model}

We model the system as a continuous-time Markov chain. The states are denoted by $(a_1,a_2,a_q,b_2,b_q)$, where $a_1$, $a_2$ and $a_q$ are the number of level 1 calls being served by level 1 agents, being served by level 2 agents, and waiting in the queue, respectively, and $b_2$ and $b_q$ are the number of level 2 calls being served by level 2 agents, and waiting in the queue, respectively.

Depending on the state, seven possible events can occur - an arrival of each level of caller, an abandonment by each level of caller, and a service of each level of caller by an agent at the same level, or a service of a level 1 caller by a level 2 agent. Whenever these events fall into one of the categories described in Subsection \ref{DM}, we need to make a decision.

For $i \in \{1,2\}$, let $k_i$ be the total number of agents at level $i$. Therefore, when the state is $\textbf{s} = (a_1,a_2,a_q,b_2,b_q)$ there are $L_1(\mathbf{s}) = k_1 - a_1$ agents available at level 1 and $L_2(\mathbf{s}) = k_2 - a_2 - b_2$ agents available at level 2. Let $\ell$ be the total capacity of the system, implying that the total number of callers at any point in the system can be at most $\ell$, including callers from all levels waiting in the queue and being served.

The state space $S_g \coloneqq \{(a_1,a_2,a_q,b_2,b_q)\}$ is constrained such that
\begin{itemize}
    \item $a_1 \in \{0,1,\dots,k_1\}$,
    \item $a_2 \in \{0,1,\dots,k_2\}$,
    \item $b_2 \in \{0,1,\dots,k_2-a_2\}$,
    \item $a_1 + a_2 + a_q + b_2 + b_q \leq \ell$,
    \item if $a_1 < k_1$, then $a_q = 0$, and
    \item if $a_2 + b_2 < k_2$, then $a_q = 0$ or $b_q = 0$.
\end{itemize}
The second last constraint rules out the case where level 1 calls are in the queue and level 1 agents are available. The last constraint rules out the case where both level 1 and level 2 calls are in the queue while level 2 agents are available.

For all $\mathbf{s} \in S_g$, let $\mathcal{E}(\mathbf{s})$ denote the set of possible events that can occur when the system is in state $\mathbf{s}$. Let $e_{ai}$ denote the arrival of a level $i$ caller, $e_{si}$ denote the completion of a service by a level $i$ agent, and $e_{abi}$ denote an abandonment by a level $i$ caller. We use $\mathcal{A}(\mathbf{s},e)$ to denote the set of possible actions that can be taken when the system is in state $\mathbf{s}$ and an event $e$ occurs. Let $A_i$ denote the action that an incoming call of level $i$ is assigned to a level $i$ agent, and $A_q$ represent the action that an incoming call is added to the queue. 
Similarly, when there is a service completion by an agent at level $i$, let $S_{i}$ represent the action that a level $i$ call from the queue is assigned to that agent, and $S_{q}$ represent the action of leaving the agent unoccupied.

For many states and events, it may happen that the action set contains only one action and does not require any decision making. For example, if the system is in a state such that $L_1(\mathbf{s})>0$, then $\mathcal{A}(\mathbf{s},e_{a1}) $ contains only $A_1$, since an arriving level 1 call will always be assigned to an available level 1 agent.

For $i \in \{1,2\}$, let $\lambda_i$ and $\theta_i$ be the arrival and the abandonment rate of level $i$ calls, respectively. For $i \in \{1,2\}$, let $\mu_i$ be the service rate of level $i$ calls by level $i$ agents and let $\mu_1'$ be the service rate of level 1 calls by level 2 agents. Let $\gamma_{i}$ be the cost incurred when a level $i$ call abandons. For a specific policy, which defines the action that needs to be taken in every state, it is straightforward to construct a transition matrix $Q$ for the underlying continuous-time Markov chain in terms of the parameters $\lambda_i$, $\mu_i$, $\mu_1'$ and $\theta_i$. From this, it is also straightforward to derive the transition matrix for the embedded discrete-time jump chain. An example of this construction is provided in Appendix \ref{SubS:Transition}. 

We now incorporate these policy-dependent discrete-time transition matrices into a Markov decision process formulation.

Let 
\begin{equation}
    \kappa(\mathbf{s}) = (\lambda_1 + \lambda_2)\mathbbm{1}_{a_1+a_2+a_q+b_2+b_q < \ell} + \mu_1 a_1  + \mu_1' a_2 + \mu_2 b_2 + \theta_1 a_q + \theta_2 b_q
\end{equation}
so that $\kappa(\mathbf{s})$ is the total transition rate from state $\mathbf{s}.$
For $i \in \{1,2\}$ and $\textbf{s} \in S_g$, let 
\begin{equation}\label{eq:bellman}
    p_{ai} = \lambda_i\mathbbm{1}_{a_1+a_2+a_q+b_2+b_q < \ell}/\kappa(\mathbf{s})
\end{equation} 
be the probability that a level $i$ call arrives, 
\begin{equation}
    p_{s1} = \mu_1 a_1/\kappa(\mathbf{s})
\end{equation} be the probability that a level 1 call's service is completed by a level 1 agent, 
\begin{equation}
    p_{s1}' = \mu_1' a_2/\kappa(\mathbf{s})
\end{equation} 
be the probability that a level 1 call's service is completed by a level 2 agent, 
\begin{equation}
    p_{s2} = \mu_2 b_2/\kappa(\mathbf{s})
\end{equation} 
be the probability that a level 2 call's service is completed by a level 2 agent, 
\begin{equation}
    p_{ab1} = \theta_1 a_q/\kappa(\mathbf{s})
\end{equation} 
be the probability that a level 1 call abandons, and 
\begin{equation}\label{eq:bellman_end}
    p_{ab2} = \theta_2 b_q/\kappa(\mathbf{s})
\end{equation} 
be the probability that a level 2 call abandons.

\subsection{A discrete time approximation}\label{dta}

In this section we apply the dynamic programming method of Richard Bellman (\textcite{bellman1957dynamic}) to the discrete-time embedded chain defined by equations \eqref{eq:bellman} -- \eqref{eq:bellman_end} over a fixed total number $M$ of time steps. We use the parameter values given in Table \ref{GenEg1}, along with $\ell = 20$ and $M = 100$.
\begin{table}[h]
    \centering
    \begin{tabular}{|c|c|c|}
    \hline
         $i$ & 1 & 2 \\ \hline
         $\lambda_i$ & 1 & 1 \\
         $\mu_i$ & 1/4 & 1/4 \\
         $\mu_i'$ & 1/4 & - \\
         $\theta_i$ & 1 & 2 \\
         $k_i$ & 5 & 5 \\
         $\gamma_{i}$ & 1 & 2 \\ \hline
    \end{tabular}
    \caption{Parameter values for the discrete-time example.}
    \label{GenEg1}
\end{table}  

For $m \in \{0,1,\dots,M\}$ and $\mathbf{s} \in S_g$, let $V(m,\textbf{s})$ be the minimum expected cost of abandonments between time steps $m$ and $M$ given that we are in state $\mathbf{s}$ at time step $m$. Since we only want to minimise the expected cost of abandonments \textit{between} time steps $m$ and $M$, we take $V(M,\textbf{s}) = 0$ $\forall \mathbf{s} \in S_g$.

Depending on the state $\textbf{s}$ at time step $m$, there will be certain actions that need to be taken whenever there is an event. 

Consider the state $\textbf{s}= (5,1,1,2,0)$. There are five level 1 callers being served by level 1 agents, one level 1 caller being served by a level 2 agent, one level 1 caller waiting in the queue, two level 2 callers being served by level 2 agents, two level 2 agents unoccupied and no level 2 callers waiting in the queue. In this state, there are two possible events that require a decision -- the arrival of a level 1 caller or the completion of a service by a level 2 agent. For each of these events $e$, there are two possible actions in the set $\mathcal{A}(\mathbf{s},e)$, resulting in a total of four possible sets of actions.

\begin{itemize}
    \item The first set of actions is to assign an arriving level 1 caller to an available level 2 agent, and to assign a level 1 caller waiting in the queue to a level 2 agent if there is a service completion by a level 2 agent. Under this set of actions,
    \begin{equation}\footnotesize
    \label{V1}
        \begin{split}
            V(m,(5,1,1,2,0)) &= p_{a1} V(m+1,(5,2,1,2,0)) + p_{a2} V(m+1,(5,1,1,3,0)) +\\& p_{s1} V(m+1,(5,1,0,2,0)) + p_{s1}' V(m+1,(5,1,0,2,0)) +\\& p_{s2} V(m+1,(5,2,0,1,0)) + p_{ab1} (V(m+1,(5,1,0,2,0))+\gamma_{1}),
        \end{split}
    \end{equation}
    
    \item The second set of actions is to send an arriving level 1 caller to the queue, and to assign a level 1 caller waiting in the queue to a level 2 agent if there is a service completion by a level 2 agent. Under this set of actions,
    \begin{equation}\footnotesize
    \label{V2}
        \begin{split}
             V(m,(5,1,1,2,0)) &= p_{a1} V(m+1,(5,1,2,2,0)) + p_{a2} V(m+1,(5,1,1,3,0)) +\\& p_{s1} V(m+1,(5,1,0,2,0)) + p_{s1}' V(m+1,(5,1,0,2,0)) +\\& p_{s2} V(m+1,(5,2,0,1,0)) + p_{ab1} (V(m+1,(5,1,0,2,0))+\gamma_{1}),
        \end{split}
    \end{equation}
    \item The third set of actions is to assign an arriving level 1 caller to an available level 2 agent, but not to assign a level 1 caller waiting in the queue to a level 2 agent when there is a service completion by a level 2 agent. Under this set of actions,
    \begin{equation}\footnotesize
    \label{V3}
        \begin{split}
             V(m,(5,1,1,2,0)) &= p_{a1} V(m+1,(5,2,1,2,0)) + p_{a2} V(m+1,(5,1,1,3,0)) +\\& p_{s1} V(m+1,(5,1,0,2,0)) + p_{s1}' V(m+1,(5,0,1,2,0)) +\\& p_{s2} V(m+1,(5,1,1,1,0)) + p_{ab1} (V(m+1,(5,1,0,2,0))+\gamma_{1}),
        \end{split}
    \end{equation}
    \item The fourth set of actions is to send an arriving level 1 caller to the queue, and not to assign a level 1 caller waiting in the queue to a level 2 agent when there is a service completion by a level 2 agent. Under this set of actions,
    \begin{equation}\footnotesize
    \label{V4}
        \begin{split}
             V(m,(5,1,1,2,0)) &= p_{a1} V(m+1,(5,1,2,2,0)) + p_{a2} V(m+1,(5,1,1,3,0)) +\\& p_{s1} V(m+1,(5,1,0,2,0)) + p_{s1}' V(m+1,(5,0,1,2,0)) +\\& p_{s2} V(m+1,(5,1,1,1,0)) + p_{ab1} (V(m+1,(5,1,0,2,0))+\gamma_{1}).
        \end{split}
    \end{equation}
    
\end{itemize}

For given values of $V(m+1,\textbf{s})$ on the right hand sides of equations \eqref{V1} -- \eqref{V4}, we select the set of actions that results in the minimum value of $V(m,(5,1,1,2,0))$. By proceeding backwards in this way from time step $M-1$ to $0$, and considering all states, we obtain the optimal action for each state and time step. Together, these optimal actions form the optimal policy. 
This policy can be described in an action matrix with the number of rows  equal to the size of the state space and the number of columns equal to $M+1$.

For the set of parameters given in Table \ref{GenEg1}, we first discuss the optimal policy at $m = 0$. We find that it is optimal to assign level 2 calls to level 2 agents whenever there is a level 2 agent available irrespective of the number of level 1 callers in the queue.

When all agents are busy and both level 1 and level 2 calls are waiting in the queue, it is optimal to assign a level 2 call to a level 2 agent who becomes free, and when only level 1 calls are in the queue, the optimal action is to leave all the level 1 calls in the queue level 1 irrespective of how many there are.

When all level 1 agents are busy and level 2 agents are available, there are two subcases -- when no one is waiting in the queue, we assign an incoming level 1 call to a level 2 agent only if there are more than two level 2 agents available, and when level 1 calls are waiting in the queue, the optimal policy is described in Table \ref{action1} ($L_2(\mathbf{s})$ is the number of level 2 agents available out of 5 and $a_q$ is the number of level 1 callers in the queue).  
\begin{table}[h]
    \centering
    \begin{tabular}{|p{1cm}|p{1cm}|p{9cm}|}
    \hline
    $L_2(\mathbf{s})$ & $a_q$ & Optimal action\\\hline
    3-5 & 1 & Assign a level 1 call to a level 2 agent when there is an arrival or a service.\\
    2 & 1 & Do not assign a level 1 call to a level 2 agent when there is an arrival but assign a level 1 call to a level 2 agent when there is a service.\\
    1 & 1 & Do not assign a level 1 call to a level 2 agent in the case of an arrival or service.\\ \hline
    2-5 & 2 & Assign a level 1 call to a level 2 agent when there is an arrival or a service.\\
    1 & 2 & Do not assign a level 1 call to a level 2 agent in the case of an arrival or service.\\ \hline
    2-5 & 3-10 & Assign a level 1 call to a level 2 agent when there is an arrival or a service.\\
    1 & 3-10 & Do not assign level 1 call to a level 2 agent when there is an arrival but assign a level 1 call to a level 2 agent when there is a service.\\ \hline
    1-4 & 11-15 & Assign a level 1 call to a level 2 agent when there is an arrival or a service.\\ \hline
    \end{tabular}
    \caption{Optimal actions for the case when all the level 1 agents are busy, some level 2 agents are available and some level 1 calls are waiting in the queue.}
    \label{action1}
\end{table}

As can be observed, the policy becomes more and more flexible (or requires less reservation of level 2 servers) as the number of level 1 calls in the queue increases. When there is only one level 1 customer in the queue, we start with reserving two agents at level 2 for level 2 calls only, then we decrease the reservation to one and finally, we remove it completely. Note that in the last row of the table, although we have written $11-15$ for $a_q$, the number of callers in the queue can increase only to the level allowed by the maximum capacity. For instance, if all level 1 agents and three level 2 agents are busy, this implies that there are already $8$ callers in the system, so the number of callers in the queue can increase only to $12$.

As we get closer to the time horizon, the optimal policy becomes more flexible in allocating level 1 calls to level 2 servers. The policy remains as in Table \ref{action1} for the first $72$ steps, and then it changes slightly. After $87$ steps, it starts assigning level 1 calls to level 2 agents whenever there is more than one level 2 agent available, and for the last three time steps, it assigns level 1 calls to level 2 agents whenever there is any availability. 

As we can see in Table \ref{GenEg1}, the abandonment rate and the cost of abandonment for level 2 calls are both twice that of level 1 calls. If we consider examples where the difference in abandonment rates and costs between the two levels is not as significant, the optimal policy obtained from this method is to assign a level 1 call to a level 2 agent whenever a level 2 agent is available, and to prioritise assigning a level 2 call when both types of calls are present. Later in the paper, we will refer to such a policy as the {\it zero-reservation policy}.

The overall pattern of the optimal policy that we observe is that if the number of level 1 calls waiting in the queue is small, then we should reserve more level 2 agents exclusively for level 2 calls. As the number of level 1 calls waiting increases, we should assign level 1 calls to level 2 agents with more flexibility. The optimal policy also depends on how close the time horizon: closer to the time horizon, the optimal policy becomes increasingly flexible in assigning level 1 callers to level 2 agents.

\subsection{An infinite time horizon Markov decision process with discounting}\label{cta}

In this section, in order to get some more insight into the type of policies that are optimal for our system, we find the optimal stationary policy for an approximate infinite horizon continuous-time Markov decision process with discounting  (\textcite{puterman2014markov}, \textcite{xie2019optimizing}, \textcite{baird1994reinforcement}).

Since we are now working with an infinite time horizon, the value function depends only on the state and not on the time. For a positive discounting factor $\delta$, let $C_\delta(\mathbf{s})$ be a random variable denoting the discounted number of abandonments over the infinite horizon given that we start in state $\mathbf{s} \in S_g$  and $V(\mathbf{s}) = \mathbbm{E}(C_\delta(\mathbf{s}))$. Let $T_0$ be the time to the next event, then $T_0$ is a random variable distributed exponentially with rate $\kappa(\mathbf{s})$. Now, we can write expressions for the expected cost of abandonments given the time to the next event $\mathbbm{E}(C_\mathbf{s}|T_0 = t)$ to calculate $V(\mathbf{s})$.

For example, again consider state $(5,1,1,2,0)$. We saw in Section \ref{dta} that, for this state, there are two possible actions when a level 1 call arrives, and two actions possible when a level 2 agent completes a service. Hence,
\begin{equation}
\begin{split}    
    \mathbbm{E}(C_\delta(5,1,1,2,0)|T_0 = t) =& e^{-\delta t}(p_{a1} \min(V(5,2,1,2,0),V(5,1,2,2,0)) +\\& p_{a2} V(5,1,1,3,0) + p_{s1} V(4,1,1,2,0) +\\& \min(p_{s1}' V(5,1,0,2,0) + p_{s2} V(5,2,0,1,0),\\&p_{s1}' V(5,0,1,2,0) + p_{s2} V(5,1,1,1,0)) +\\& p_{ab1} (V(5,1,0,2,0)+\gamma_1)).
\end{split}
\end{equation}

Now, using the independence of $T_0$ and the subsequent evolution of the process, we can write
\begin{equation}\label{eq:disc}
\begin{split}    
    V(5,1,1,2,0) =& \left(\frac{\kappa(5,1,1,2,0)}{\delta + \kappa(5,1,1,2,0)}\right)(p_{a1} \min(V(5,2,1,2,0),V(5,1,2,2,0)) \\& + p_{a2} V(5,1,1,3,0) + p_{s1} V(5,1,0,2,0) \\& +\min(p_{s1}' V(5,1,0,2,0) + p_{s2} V(5,2,0,1,0),\\&p_{s1}' V(5,0,1,2,0) + p_{s2} V(5,1,1,1,0)) \\& + p_{ab1} (V(5,1,0,2,0)+\gamma_1)).
\end{split}
\end{equation}
Similarly we can write equations for the value function in every state. These equations are non-linear because they contain $\min(\cdot)$. We use an iterative method to obtain an approximate solution (\textcite{powell2007approximate}, \textcite[Chapter 1]{chang2013simulation}). 

Let $K$ be a fixed number of iterations and $\alpha \in [0,1]$ be a step size. Then, for $k \in \{0,1,\dots,K-1\}$, we write 
\begin{equation}\label{eq:iter}\footnotesize
\begin{split}    
    V_{k+1}(5,1,1,2,0) =& \alpha V_{k}(5,1,1,2,0)\\& +(1-\alpha)\left[\left(\frac{\kappa(5,1,1,2,0)}{\delta +  \kappa(5,1,1,2,0)}\right)(p_{a1} \min(V_k(5,2,1,2,0),V_k(5,1,2,2,0))\right. \\& + p_{a2} V_k(5,1,1,3,0) + p_{s1} V_k(5,1,0,2,0)\\& +\min(p_{s1}' V_k(5,1,0,2,0) + p_{s2} V_k(5,2,0,1,0),\\&p_{s1}' V_k(5,0,1,2,0) + p_{s2} V_k(5,1,1,1,0)) \\& \left. + p_{ab1} (V_k(5,1,0,2,0)+\gamma_1))\right],
\end{split}
\end{equation}
with analogous equations for every other state.
We can initialise the iteration by choosing the $V_0(\mathbf{s})$ arbitrarily. The values of $\alpha$ and $V_0(\mathbf{s})$ will affect the convergence speed of the algorithm, but not the final values $V_K(\mathbf{s})$ and the optimal policy. We can also increase the computational efficiency of this method by introducing a stopping criterion that compares the vector of value functions at iteration $k+1$ with that at iteration $k$. We obtain the optimal stationary policy by extracting the arguments that lead to the minimum values on the right hand side.

\subsubsection{Example}

We tested the algorithm with two different initialisations of $V_0(\mathbf{s})$ -- one using the values from the model discussed in Section \ref{dta} and the other just simply putting them equal to zero. We also tested it with different values of $\alpha$. All combinations led to the same optimal policy and very close value functions with $K = 1000$. 

We used the above procedure to calculate the optimal policy, again using the parameter values in Table \ref{GenEg1}. Again, it is always optimal to assign level 2 calls to level 2 agents whenever there is an availability. When all the level 1 agents are busy while level 2 agents are available and no one is waiting in the queue, we assign an incoming level 1 call to a level 2 agent only if there are more than two level 2 agents available. In the case when level 1 calls are waiting in the queue, the policy is described in Table \ref{action2}.

\begin{table}[h]
    \centering
    \begin{tabular}{|p{1cm}|p{1cm}|p{9cm}|}
    \hline
    $L_2(\mathbf{s})$ & $a_q$ & Optimal action\\\hline
    2-5 & 1 & Assign a level 1 call to a level 2 agent when there is an arrival or a service.\\
    1 & 1 & Do not assign a level 1 call to a level 2 agent in the case of an arrival or service.\\ \hline
    2-5 & 2-10 & Assign a level 1 call to a level 2 agent when there is an arrival or a service.\\
    1 & 2-10 & Do not assign level 1 call to a level 2 agent when there is an arrival but assign a level 1 call to a level 2 agent when there is a service.\\ \hline
    1-5 & 11-15 & Assign a level 1 call to a level 2 agent when there is an arrival or a service.\\ \hline
    \end{tabular}
    \caption{Optimal actions for the case when all the level 1 agents are busy, some level 2 agents are available, and some level 1 calls are waiting in the queue.}
    \label{action2}
\end{table}

When all agents are busy and both level 1 and level 2 calls are waiting in the queue, it is optimal to assign a level 2 call to a level 2 agent who becomes free. When only level 1 calls are in the queue, the optimal action is not to assign a level 1 call to a level 2 agent from the queue regardless of how many are waiting. We observe that the results are similar to those obtained from the discrete time model in Table \ref{action1}.

\section{Reservation policies}\label{ch:res}

We can characterise the optimal policies that we derived in Sections \ref{dta} and \ref{cta} as {\it bi-threshold policies} in which we take actions depending on the number of free agents available at level 2 and the number of level 1 calls waiting in the queue. These policies assign level 1 calls to level 2 servers if both the numbers of unoccupied level 2 agents and level 1 calls waiting in the queue are above respective thresholds. 

We also observe from the finite-horizon model in Section \ref{dta} that in finite time the optimal policy changes in favour of allocating level 1 calls to level 2 servers as the number of time steps left until the planning horizon decreases. 

It is reasonable to say that such a policy could be too complicated to implement in practice. It is computationally expensive as it requires a large number of state space variables, even for a two-level problem  and the exact policy is sensitive to changes in the parameter values. 

Based upon the intuition gained from looking at the models in Sections \ref{dta} and \ref{cta}, we propose a class of policies that are likely to be close to be optimal for hierarchical models with $\nu \geq 2$. We call these policies \textit{reservation policies}. They will form the focus of the rest of this paper.

\begin{definition}\label{Def}
A \textit{reservation policy} for a hierarchical model of the class defined in Section \ref{intro} is a policy under which calls are assigned to agents based on reservation thresholds $\Theta_i$ for each level $i$ agent. Under this policy

\begin{itemize}
\item if there is an arrival of a level $i-1$ caller,
\begin{itemize}
\item first, attempt to assign the caller to an available level $i-1$ agent,
\item if all level $i-1$ agents are busy, then check the number of available level $i$ agents,
\item if at least $\Theta_i$ level $i$ agents are available, assign the call to one of them,
\item if fewer than $\Theta_i$ level $i$ agents are available, assign the caller to the queue.
\end{itemize}
\item if a level $i$ agent completes a service,
\begin{itemize}
    \item the agent first checks for level $i$ callers in the queue and serves such a caller if any exist
    \item if there are no level $i$ callers in the queue, but there are level $i-1$ callers waiting -- the agent will serve a level $i-1$ caller if the number of available level $i$ agents (including the newly available agent) is at least $\Theta_i$.
\end{itemize}
\end{itemize}
\end{definition}

We see that a reservation policy depends only on the number of higher level agents available and not on the number of lower level callers callers in the queue. Such a policy is parameterised by the values $\mathbf{\Theta} = (\Theta_2,\ldots,\Theta_\nu)$ of the thresholds.

Using the method that we develop in Section \ref{ab4} for a model with $\nu=2$, we analytically calculated the expected cost of abandonments $C(\Theta_2,t)$ for the parameter values in Table \ref{GenEg1} and time duration $(0,t)$ for three different values of $t$ and all the possible values of $\Theta_2 \in \{0,1,\dots,k_2\}$. The values are given in Table \ref{lv2eg}.

\begin{table}[h!]
\centering
\begin{subtable}{0.3\textwidth}
\centering
\begin{tabular}{|c|c|}
\hline
$\Theta_2$ & $C(\Theta_2,2)$ \\
\hline
0 & 0.0078 \\
1 & 0.0078 \\
2 & 0.0081 \\
3 & 0.0086 \\
4 & 0.0094 \\
5 & 0.0100 \\
\hline
\end{tabular}
\caption{$t = 2$}
\end{subtable}
\begin{subtable}{0.3\textwidth}
\centering
\begin{tabular}{|c|c|}
\hline
$\Theta_2$ & $C(\Theta_2,15)$ \\
\hline
0 & 4.467 \\
1 & 4.278 \\
2 & 4.290 \\
3 & 4.387 \\
4 & 4.476 \\
5 & 4.511 \\
\hline
\end{tabular}
\caption{$t = 15$}
\end{subtable}
\begin{subtable}{0.3\textwidth}
\centering
\begin{tabular}{|c|c|}
\hline
$\Theta_2$ & $C(\Theta_2,60)$ \\
\hline
0 & 27.071 \\
1 & 25.706 \\
2 & 25.528 \\
3 & 25.864 \\
4 & 26.204 \\
5 & 26.332 \\
\hline
\end{tabular}
\caption{$t = 60$}
\end{subtable}
\caption{Expected cost of abandonments for parameter values given in Table \ref{GenEg1}.}
    \label{lv2eg}
\end{table}

We observe that for small values of $t$ (analogous to being near the finite horizon in the discrete time method), the minimum expected cost of abandonments occurs at $\Theta_2 = 0$, that is, the optimal policy is the zero reservation policy where we do not reserve any level 2 agents and assign a level 1 call to a level 2 agent whenever there is an availability. For $t=15$ and values around that, the minimum is achieved at $\Theta_2 = 1$ implying that it is optimal to reserve one level 2 agent for level 2 calls only, and for the larger values of $t$, it is optimal to reserve two agents at level $2$. This is similar to what we observe in Section \ref{dta}. 

Overall, when we compare the optimal reservation policy with the optimal policies obtained in Sections \ref{dta} and \ref{cta}, we find that the differences arise primarily in scenarios where a large number of level 1 callers are waiting in the queue.

We implemented the policies derived in Section \ref{dta} at $m = 0$ and Section \ref{cta}, along with the reservation policy with $\Theta_2 = 2$ on a simulated underlying continuous-time Markov chain using the parameter values in Table \ref{GenEg1}. We observed that the expected cost of abandonments is nearly identical across all three policies and the performance was not distinguishable just on the basis of the simulation results.

Moreover, when the abandonment rates and the associated costs are similar across the two levels, the optimal reservation policy corresponds to a zero reservation policy, the same as the policy derived in Sections \ref{dta} and \ref{cta}.

In Section \ref{ab4}, we shall present a method to calculate the exact expected cost of abandonments for a reservation policy in continuous time for a finite time horizon. Since reservation policies require a smaller number of variables compared to the policies discussed in Section \ref{General}, we are able to extend our analysis to the case where there are four levels of both agents and callers. 

\section{Reservation policies for call centres with four levels} \label{ab4}

In this section, we will formulate a modification of the method introduced in \textcite{chiera_what_2002} for our industry partner's call centre which has four levels of calls and agents.
The objective will be to calculate the optimal thresholds $\mathbf{\Theta}$ and the corresponding expected costs of abandonments for the class of reservation policies defined in \ref{Def}.

We model the system as a continuous-time Markov chain, assuming that the calls are arriving according to Poisson processes, and the service and abandonment times are exponentially distributed.

For $i \in \{1,2,3,4\}$, let $\lambda_i$, $\theta_i$, and $\gamma_i$ be the arrival rate, abandonment rate, and abandonment cost of level $i$ callers, respectively. Let $k_i$ be the total number of level $i$ agents, $\mu_i$ be the service rate when level $i$ callers are served by level $i$ agents, and for $i \in \{1,2,3\}$, let $\mu_i'$ be the service rate when level $i$ callers are served by level $i+1$ agents. Let $\ell$ be the total capacity of the system shared by all callers regardless of their level and let $\beta$ be the cost incurred when customer is lost/blocked. Because the number of reserved agents has to be less than the total number of agents at that level, we have for $i \in \{2,3,4\}$, $\Theta_i \in \{0,1,\dots,k_i\}$.

\subsection{State space}

We denote the state of the system by $(a,a_1,b,b_1,c,c_1,d)$, where $a, b, c$ and $d$ are the number of level 1, 2, 3 and 4 callers in the system, and $a_1, b_1$ and $c_1$ are the number of level 1, 2 and 3 callers in the system being served by level 2, 3 and 4 agents, respectively. This state space description is slightly different from that used in Section \ref{gen_model}, because the nature of the reservation policy allows us to convey the same amount of information using fewer variables.

When the system is in state $(a,a_1,b,b_1,c,c_1,d)$ there are $a$ level 1 callers, $a_1$ of which are being served by level 2 agents. Hence, $\min(a-a_1,k_1)$ level 1 callers are being served by level 1 agents, and $\max(0,a-(a_1+k_1))$ level 1 callers are waiting in the queue. It is possible that some of the level 1 agents are free while level 2 agents are serving level 1 callers if level 1 agents became free after the calls were already assigned to level 2 agents. Now, if $a_1$ level 2 agents are serving level 1 callers, and $b_1$ level 2 callers are being served by level 3 agents, then this implies $\min(b-b_1,k_2-a_1)$ level 2 callers are being served by level 2 agents, and $\max(0,b-(b_1+k_2-a_1))$ are waiting in the queue. Similarly, we can calculate how many level 3 and level 4 callers are waiting and how many are being served by each level of agent.

Here, $a,b,c,d \in \{0,1,\dots,\ell\}$ and $a+b+c+d \leq \ell$. The maximum value $a_1$ can take is $\min(k_2-\Theta_2,a)$ because by the definition of $a_1$, it cannot be more than the total number of level 1 callers and has to be less than or equal to the total number of level 2 agents. The minimum value of $a_1$ is zero in all cases except for when $a > k_1$ and $b < k_2-\Theta_2$ (that is, the number of level 1 callers are more than the number of level 1 agents and there are more than $\Theta_2$ level 2 agents available). Hence,\begin{equation}\min(a-k_1,k_2-b-\Theta_2) \leq a_1 \leq \min(k_2-\Theta_2,a).\end{equation} Similarly, \begin{equation}\min(b-k_2,k_3-c-\Theta_3)\mathbbm{1}_{b>k_2,c<k_3-\Theta_3} \leq b_1 \leq \min(k_3-\Theta_3,b)\end{equation} and \begin{equation}\min(c-k_3,k_4-d-\Theta_4)\mathbbm{1}_{c>k_3,d<k_4-\Theta_4} \leq c_1 \leq \min(k_4-\Theta_4,c).\end{equation}
Hence, the state space is given by $S = \{(a,a_1,b,b_1,c,c_1,d)\}$ such that
\begin{itemize}
    \item $a,b,c,d \in \{0,1,\dots,l\}$,
    \item $a+b+c+d \leq \ell$,
    \item $\min(a-k_1,k_2-b-\Theta_2) \mathbbm{1}_{a > k_1,b < k_2-\Theta_2} \leq a_1 \leq \min(k_2-\Theta_2,a)$, 
    \item $\min(b-k_2,k_3-c-\Theta_3)\mathbbm{1}_{b>k_2,c<k_3-\Theta_3} \leq b_1 \leq \min(k_3-\Theta_3,b)$ and 
    \item $\min(c-k_3,k_4-d-\Theta_4)\mathbbm{1}_{c>k_3,d<k_4-\Theta_4} \leq c_1 \leq \min(k_4-\Theta_4,c)$.
\end{itemize}
The details of the transition rates for this model are given in Appendix \ref{SubS:Transition}.

\subsection{The exact expected cost of losses and abandonments}

For a given state $(a,a_1,b,b_1,c,c_1,d) \in S$ at time $0$ and reservation vector $\mathbf{\Theta}$, define $C_{a,a_1,b,b_1,c,c_1,d}(\mathbf{\Theta},t), t \geq 0$ to be the expected cost of abandonments and losses during the time interval $(0,t)$.

Using a methodology similar to that used in \textcite{chiera_what_2002}, we condition on the first time that the system leaves its original state, writing an expression for $C_{a,a_1,b,b_1,c,c_1,d}(\mathbf{\Theta},t|x)$ where $x$ is the sojourn time in the state  $(a,a_1,b,b_1,c,c_1,d)$. 

If $t < x$, the system stays in its original state for the whole interval and there are no abandonments in time $(0,t)$. However, if the system is at its capacity, there might be lost customers. The expected cost in that case is given by $ t\, \beta\,\sum_{i = 1}^4 \lambda_i$. 

If $t\geq x$ and, for example, the system is currently in a state such that, $a < k_1, b+a_1 < k_2, c+b_1 < k_3, d+c_1 > k_4$, implying that there are agents available at each level except level 4, and there are customers waiting at level 4, then the expected cost of abandonments and losses is given by 

\begin{equation}
    C_{a,a_1,b,b_1,c,c_1,d}(\mathbf{\Theta},t|x) = \frac{C^n_{a,a_1,b,b_1,c,c_1,d}(\mathbf{\Theta},t|x)}{\kappa_{a,a_1,b,b_1,c,c_1,d}(\mathbf{\Theta})}
\end{equation}
where
\begin{equation}
\label{eq:cond}
\begin{split}    
C^n_{a,a_1,b,b_1,c,c_1,d}(\mathbf{\Theta},t|x) = & \lambda_1 C_{a+1,a_1,b,b_1,c,c_1,d}(\mathbf{\Theta},t-x)\\& + \lambda_2 C_{a,a_1,b+1,b_1,c,c_1,d}(\mathbf{\Theta},t-x) \\& + \lambda_3 C_{a,a_1,b,b_1,c+1,c_1,d}(\mathbf{\Theta},t-x) \\& + \lambda_4 C_{a,a_1,b,b_1,c,c_1,d+1}(\mathbf{\Theta},t-x) \\&
+ (a-a_1)\mu_1 C_{a-1,a_1,b,b_1,c,c_1,d}(\mathbf{\Theta},t-x)\\& + (b-b_1)\mu_2 C_{a,a_1,b-1,b_1,c,c_1,d}(\mathbf{\Theta},t-x)\\& + (c-c_1)\mu_3 C_{a,a_1,b,b_1,c-1,c_1,d}(\mathbf{\Theta},t-x)\\& + (k_4 - c_1)\mu_4 C_{a,a_1,b,b_1,c,c_1,d-1}(\mathbf{\Theta},t-x)\\& 
+ a_1 \mu_1' C_{a-1,a_1-1,b,b_1,c,c_1,d}(\mathbf{\Theta},t-x)\\& + b_1 \mu_2' C_{a,a_1,b-1,b_1-1,c,c_1,d}(\mathbf{\Theta},t-x)\\& + c_1 \mu_3' C_{a,a_1,b,b_1,c-1,c_1-1,d}(\mathbf{\Theta},t-x)\\&
+ (d-(k_4-c_1))\theta_4(\gamma_4+C_{a,a_1,b,b_1,c,c_1,d-1}(\mathbf{\Theta},t-x))
\end{split}
\end{equation}
and
\begin{equation}
\begin{split}
    \kappa_{a,a_1,b,b_1,c,c_1,d}(\mathbf{\Theta}) = & \sum_{i = 1}^4\lambda_i + (a-a_1)\mu_1 + (b-b_1)\mu_2 + (c-c_1)\mu_3 + (k_4-c_1)\mu_4 \\& + a_1\mu_1' + b_1\mu_2' + c_1\mu_3' + (d-(k_4-c_1))\theta_4.
\end{split}
\end{equation}
If the system is at the capacity this equation needs to be modified to include the blocking cost in time $(0,x)$. After writing similar expressions for all states $C_{a,a_1,b,b_1,c,c_1,d}(\mathbf{\Theta},t|x)$, we can integrate to remove the conditioning on $x$, observing that $x$ is an exponentially distributed random variable. The right hand side of \eqref{eq:cond} becomes a sum of convolution integrals which suggests that we should take Laplace transforms. This results in expressions of the form
\begin{equation}
\label{eq:LT}
\Tilde{C}_{a,a_1,b,b_1,c,c_1,d}(\mathbf{\Theta},s) = \dfrac{\xi_{a,a_1,b,b_1,c,c_1,d}(\mathbf{\Theta},s)}{\psi_{a,a_1,b,b_1,c,c_1,d}(s)}
\end{equation}
for the Laplace transform of the expected cost of losses and abandonments in time $(0,t)$, conditional on the system being in state $(a,a_1,b,b_1,c,c_1,d)$ at time zero. 

The exact expressions for $\xi_{a,a_1,b,b_1,c,c_1,d}(\mathbf{\Theta},s)$ and $\psi_{a,a_1,b,b_1,c,c_1,d}(s)$
are given in Appendix \ref{laplace}. 

We can solve equations \eqref{eq:LT} numerically and then invert them using the Euler method discussed in \textcite{abate1995numerical}.

A similar approach to that described above was used by \textcite{knessl2015transient} for a transient analysis of the $M/M/m + M$ queue. However, the reservation policy component in our model makes the Laplace transform equations that we obtain more complicated.

If we let $\beta = 0$, $\gamma_i = 1$, and $\gamma_j = 0$ for $i \neq j$, we get an expression for the expected number of abandonments by level $i$ callers and if we let $\gamma_i = 0$ and $\beta = 1$ we get an expression for the expected number of lost customers.

Using the same logic as above, we can also calculate the expected waiting time of customers, which is proportional to the expected cost of abandonments with $\beta = 0$. This is not surprising. While class $i$ customers are waiting, they each have their own `exponential abandonment process’ going, which has rate $\theta_i$ and is independent of what other customers are doing. We can also calculate the expected number of services and the expected waiting time of the $j^{\text{th}}$ customer in the queue.

\subsection{Numerical results}\label{num}

In this section, we will look at some examples on how to decide the best allocation policy.

\subsubsection{Example 1}\label{eg1}

Consider a model with parameter values given in Table \ref{eg1t} along with $\beta = 0$, $\ell = 10$ and whose initial state is $(0,0,0,0,0,0,0)$.

\begin{table}[h]
    \centering
    \begin{tabular}{|c|c|c|c|c|}\hline
        $i$ & 1 & 2 & 3 & 4 \\ \hline
        $\lambda_i$ & 1 & 1/2 & 1/5 & 1/8 \\[0.2ex]
        $\mu_i$ & 2/3 & 1/2 & 1/4 & 1/6 \\ [0.2ex]
        $\mu_i'$ & 2/3 & 1/2 & 1/4 & - \\
        $\theta_i$ & 2 & 1 & 1 & 1 \\
        $k_i$ & 3 & 2 & 2 & 2\\
        $\gamma_i$ & 1 & 1 & 1 & 1\\
        \hline
    \end{tabular}
    \caption{Parameter values for Example 1.}
    \label{eg1t}
\end{table}

\begin{table}[h!]
\centering
\begin{subtable}{0.45\textwidth}
\begin{tabular}{|c|c|c|c|}
\hline
$\Theta_2$ & $\Theta_3$ & $\Theta_3$ & $C(\mathbf{\Theta},60)$ \\
\hline
0 & 0 & 0 & 3.53 \\
0 & 0 & 1 & 4.03 \\
0 & 0 & 2 & 4.81 \\
0 & 1 & 0 & 5.22 \\
0 & 1 & 1 & 5.62 \\
0 & 1 & 2 & 6.21 \\
0 & 2 & 0 & 7.23 \\
0 & 2 & 1 & 7.57 \\
0 & 2 & 2 & 8.10 \\
1 & 0 & 0 & 5.37 \\
1 & 0 & 1 & 5.86 \\
1 & 0 & 2 & 6.59 \\
1 & 1 & 0 & 6.81 \\
1 & 1 & 1 & 7.19 \\
\hline
\end{tabular}
\end{subtable}
\hfill
\begin{subtable}{0.45\textwidth}
\begin{tabular}{|c|c|c|c|}
\hline
$\Theta_2$ & $\Theta_3$ & $\Theta_3$ & $C(\mathbf{\Theta},60)$ \\
\hline
1 & 1 & 2 & 7.77 \\
1 & 2 & 0 & 8.51 \\
1 & 2 & 1 & 8.86 \\
1 & 2 & 2 & 9.40 \\
2 & 0 & 0 & 7.13 \\
2 & 0 & 1 & 7.61 \\
2 & 0 & 2 & 8.33 \\
2 & 1 & 0 & 8.45 \\
2 & 1 & 1 & 8.84 \\
2 & 1 & 2 & 9.42 \\
2 & 2 & 0 & 10.04 \\
2 & 2 & 1 & 10.38 \\
2 & 2 & 2 & 10.92 \\
\hline
\end{tabular}
\end{subtable}
\caption{Expected cost of abandonments for the parameter values given in Table \ref{eg1t}.}
\label{t1}
\end{table}

Assuming that the unit of time is one minute, the values in Table \ref{eg1t} imply that calls are arriving one per minute, one per two minutes, one per five minutes and one per eight minutes for level $1$, $2$, $3$ and $4$ respectively, capturing the fact that more complex calls are rarer. The service rate for the four levels of calls is such that the expected service time increases with the increase in the level, and is indifferent to whether the call is being served by the same level agent or by a higher level agent. Similarly, the abandonment rate for level 1 calls is two per minute and for all other levels is one per minute. The cost of abandonment is equal for all four levels of calls and we are not counting the blocked customers (the expected number of callers blocked in this example is less than $0.3\%$ showing that $\ell = 10$ is reasonable). 

We calculated the expected number of call abandonments between $t = 0$ and $t = 60$ minutes. Table \ref{t1} gives the expected cost of abandonments for all possible values of ${\mathbf{\Theta}}$.
As can be seen from the table, the optimal value of ${\mathbf{\Theta}}$ is $(0,0,0)$, that is, not to use reservation of agents at all and allocate lower level calls to higher level agents as they are available. In this example, the optimal expected cost of abandonments is $3.53$. As a comparison, we can observe that, for the worst policy (when there is complete reservation), the expected cost of abandonments is $10.92$.

\subsubsection{Example 2}\label{eg2}

In this example, we look at a scenario in which level 4 calls are much more important than other calls, and if a level 4 caller abandons the queue it is four times more expensive than an abandonment by any other level of customer.

The parameter values are given in Table \ref{eg2t} along with $\beta = 0$ and $\ell = 10$. Table \ref{t2} gives the expected cost of abandonments between $t = 0$ and $t = 60$ minutes for all possible values of $\mathbf{\Theta}$.
\begin{table}[h]
    \centering
    \begin{tabular}{|c|c|c|c|c|}\hline
        $i$ & 1 & 2 & 3 & 4 \\ \hline
        $\lambda_i$ & 1 & 1/2 & 1/4 & 1/4 \\[0.2ex]
        $\mu_i$ & 2/3 & 1/2 & 1/4 & 1/2 \\ [0.2ex]
        $\mu_i'$ & 2/3 & 1/2 & 1/4 & - \\
        $\theta_i$ & 2 & 1 & 1 & 1 \\
        $k_i$ & 3 & 2 & 2 & 2\\
        $\gamma_i$ & 1 & 1 & 1 & 4\\
        \hline
    \end{tabular}
    \caption{Parameter values for Example 2.}
    \label{eg2t}
\end{table}

\begin{table}[h!]
\centering
\begin{subtable}{0.45\textwidth}
\centering
\begin{tabular}{|c|c|c|c|}
\hline
$\Theta_2$ & $\Theta_3$ & $\Theta_3$ & $C(\mathbf{\Theta},60)$ \\
\hline
0 & 0 & 0 & 7.78 \\
0 & 0 & 1 & 7.13 \\
0 & 0 & 2 & 7.79 \\
0 & 1 & 0 & 9.13 \\
0 & 1 & 1 & 8.63 \\
0 & 1 & 2 & 9.17 \\
0 & 2 & 0 & 10.78 \\
0 & 2 & 1 & 10.28 \\
0 & 2 & 2 & 10.78 \\
1 & 0 & 0 & 9.53 \\
1 & 0 & 1 & 8.90 \\
1 & 0 & 2 & 9.53 \\
1 & 1 & 0 & 10.67 \\
1 & 1 & 1 & 10.15 \\
\hline
\end{tabular}
\end{subtable}
\hfill
\begin{subtable}{0.45\textwidth}
\centering
\begin{tabular}{|c|c|c|c|}
\hline
$\Theta_2$ & $\Theta_3$ & $\Theta_3$ & $C(\mathbf{\Theta},60)$ \\
\hline
1 & 1 & 2 & 10.69 \\
1 & 2 & 0 & 12.07 \\
1 & 2 & 1 & 11.58 \\
1 & 2 & 2 & 12.06 \\
2 & 0 & 0 & 11.26 \\
2 & 0 & 1 & 10.63 \\
2 & 0 & 2 & 11.26 \\
2 & 1 & 0 & 12.30 \\
2 & 1 & 1 & 11.80 \\
2 & 1 & 2 & 12.32 \\
2 & 2 & 0 & 13.60 \\
2 & 2 & 1 & 13.10 \\
2 & 2 & 2 & 13.60 \\
\hline
\end{tabular}
\end{subtable}
\caption{Expected cost of abandonments for the parameter values given in Table \ref{eg2t}.}
\label{t2}
\end{table}

The expected cost of abandonments, in this case, ranges from $7.13$ to $13.60$, and the best reservation policy is $\mathbf{\Theta} = (0,0,1)$, that is, out of the two level 4 agents, we should reserve one for just level 4 callers and all other agents should be flexible.

\subsubsection{Example 3}\label{eg3}

Now, we look at an outlier example with cost of abandonment of level $3$ customers ten times that of others.

The parameter values are given in Table \ref{eg3t} along with $\beta = 0$, and $\ell = 10$. Table \ref{t3} gives the expected cost of abandonments between $t = 0$ and $t = 60$ minutes for all possible values of $\mathbf{\Theta}$.

\begin{table}[h]
    \centering
    \begin{tabular}{|c|c|c|c|c|}\hline
        $i$ & 1 & 2 & 3 & 4 \\ \hline
        $\lambda_i$ & 1/4 & 1/4 & 1/3 & 1/8 \\[0.2ex]
        $\mu_i$ & 1 & 1/4 & 1/2 & 1/6 \\ [0.2ex]
        $\mu_i'$ & 1 & 1/8 & 1/2 & - \\
        $\theta_i$ & 1 & 1 & 2 & 1 \\
        $k_i$ & 2 & 2 & 2 & 2\\
        $\gamma_i$ & 1 & 1 & 10 & 1\\
        \hline
    \end{tabular}
    \caption{Parameter values for Example 3.}
    \label{eg3t}
\end{table}

\begin{table}[h!]
\centering
\begin{subtable}{0.45\textwidth}
\centering
\begin{tabular}{|c|c|c|c|}
\hline
$\Theta_2$ & $\Theta_3$ & $\Theta_3$ & $C(\mathbf{\Theta},60)$ \\
\hline
0 & 0 & 0 & 8.25 \\
0 & 0 & 1 & 19.83 \\
0 & 0 & 2 & 33.10 \\
0 & 1 & 0 & 6.66 \\
0 & 1 & 1 & 14.91 \\
0 & 1 & 2 & 24.81 \\
0 & 2 & 0 & 6.27 \\
0 & 2 & 1 & 12.32 \\
0 & 2 & 2 & 19.81 \\
1 & 0 & 0 & 8.27 \\
1 & 0 & 1 & 19.81 \\
1 & 0 & 2 & 33.06 \\
1 & 1 & 0 & 6.70 \\
1 & 1 & 1 & 14.93 \\
\hline
\end{tabular}
\end{subtable}
\hfill
\begin{subtable}{0.45\textwidth}
\centering
\begin{tabular}{|c|c|c|c|}
\hline
$\Theta_2$ & $\Theta_3$ & $\Theta_3$ & $C(\mathbf{\Theta},60)$ \\
\hline
1 & 1 & 2 & 24.81 \\
1 & 2 & 0 & 6.31 \\
1 & 2 & 1 & 12.36 \\
1 & 2 & 2 & 19.85 \\
2 & 0 & 0 & 8.32 \\
2 & 0 & 1 & 19.86 \\
2 & 0 & 2 & 33.09 \\
2 & 1 & 0 & 6.76 \\
2 & 1 & 1 & 14.98 \\
2 & 1 & 2 & 24.86 \\
2 & 2 & 0 & 6.37 \\
2 & 2 & 1 & 12.41 \\
2 & 2 & 2 & 19.90 \\
\hline
\end{tabular}
\end{subtable}
\caption{Expected cost of abandonments for the parameter values given in Table \ref{eg3t}.}
\label{t3}
\end{table}

The expected cost of abandonments, in this case, is a lot more sensitive to the value of $\mathbf{\Theta}$. The minimum expected cost of abandonments is $6.27$ for $\mathbf{\Theta} = (0,2,0)$ and the maximum expected cost of abandonments is $33.10$ for $\mathbf{\Theta} = (0,0,2)$.

\subsubsection{Example 4 -- a multilingual call centre}\label{eg4}

We can extend the methods discussed here to analyse call centre systems that have different types of hierarchy. Consider, for example, a call centre that gives its customers an option of speaking to an agent skilled at conversing in English, Hindi or Mandarin.

We consider a scenario where English-speaking agents can only answer queries in English, while Hindi and Mandarin-speaking agents can answer calls both in their own language and English. The question arises as to how many of the Hindi and Mandarin-speaking agents we should reserve for customers asking for their specific language.

Here, we need to elaborate on the policy, particularly in the case when no English speaking agent is available. In the example in this section, we assume that if the number of Hindi speaking agents available above their threshold exceeds that of Mandarin speaking agents above their threshold, an English query is assigned to a Hindi speaking agent, and vice versa. If the number of available agents above the thresholds is equal for both languages, the English query is assigned to a Hindi speaking agent. This choice is arbitrary and can be modified to reflect service rates or other relevant factors such as cost or priority.

For notation purposes,  we refer to English speaking agents and English queries as level 1, Hindi speaking agents and queries as level 2, and Mandarin speaking agents and queries as level 3. The rest of the parameters remain the same as in previous examples, except that $\mu_2'$ and $\mu_3'$ now denote the service rates for English customers served by Hindi and Mandarin speaking agents, respectively. 

We consider the parameter values given in Table \ref{langt2} along with $\beta = 0$ and $\ell = 12$.

\begin{table}[h]
    \centering
    \begin{tabular}{|c|c|c|c|c|}\hline
        $i$ & $1$ & $2$ & $3$ \\ \hline
        $\lambda_i$ & 2 & 1/4 & 1/4 \\[0.2ex]
        $\mu_i$ & 1 & 1 & 1 \\ [0.2ex]
        $\mu_i$ & - & 1/4 & 1/4 \\[0.2ex]
        $\theta_i$ & 1 & 2 & 2 \\[0.2ex]
        $k_i$ & 3 & 2 & 2 \\[0.2ex]
        $\gamma_i$ & 1 & 1 & 1 \\
        \hline
    \end{tabular}
    \caption{Parameter values for Example 4.}
    \label{langt2}
\end{table}

The expected cost of abandonments for the parameter values given in Table $\ref{langt2}$ for $t = 120$ given that the system is empty at time $0$ is given in Table \ref{langtt2}.

\begin{table}[h!]
\centering
\begin{tabular}{|c|c|c|}
\hline
$\Theta_2$ & $\Theta_3$ & $C(\mathbf{\Theta},120)$ \\
\hline
0 & 0 & 15.01 \\
0 & 1 & 16.19 \\
0 & 2 & 20.13 \\
1 & 0 & 15.04 \\
1 & 1 & 14.77 \\
1 & 2 & 19.83 \\
2 & 0 & 20.13 \\
2 & 1 & 19.83 \\
2 & 2 & 26.54 \\
\hline
\end{tabular}
\caption{Expected cost of abandonments for the parameter values given in Table $\ref{langt2}$.}
\label{langtt2}
\end{table} 

In this case, we find that the optimal policy is to reserve one each of the Hindi and Mandarin-speaking agents.

The expected cost of abandonments for English, Hindi and Mandarin speaking customers separately is given in Table \ref{langeg2t}.

\begin{table}[h!]
\centering
\begin{subtable}{0.3\textwidth}
\centering
\begin{tabular}{|c|c|c|}
\hline
$\Theta_2$ & $\Theta_3$ & $C_1(\mathbf{\Theta},120)$ \\
\hline
0 & 0 & 3.57 \\
0 & 1 & 6.61 \\
0 & 2 & 10.84 \\
1 & 0 & 6.59 \\
1 & 1 & 11.35 \\
1 & 2 & 17.56 \\
2 & 0 & 10.84 \\
2 & 1 & 17.56 \\
2 & 2 & 25.76 \\
\hline
\end{tabular}
\end{subtable}
\hfill
\begin{subtable}{0.3\textwidth}
\centering
\begin{tabular}{|c|c|c|}
\hline
$\Theta_2$ & $\Theta_3$ & $C_2(\mathbf{\Theta},120)$ \\
\hline
0 & 0 & 6.58 \\
0 & 1 & 8.31 \\
0 & 2 & 8.90 \\
1 & 0 & 1.57 \\
1 & 1 & 1.82 \\
1 & 2 & 1.88 \\
2 & 0 & 0.39 \\
2 & 1 & 0.39 \\
2 & 2 & 0.39 \\
\hline
\end{tabular}
\end{subtable}
\hfill
\begin{subtable}{0.3\textwidth}
\centering
\begin{tabular}{|c|c|c|}
\hline
$\Theta_2$ & $\Theta_3$ & $C_3(\mathbf{\Theta},120)$ \\
\hline
0 & 0 & 4.86 \\
0 & 1 & 1.28 \\
0 & 2 & 0.39 \\
1 & 0 & 6.87 \\
1 & 1 & 1.59 \\
1 & 2 & 0.39 \\
2 & 0 & 8.90 \\
2 & 1 & 1.88 \\
2 & 2 & 0.39 \\
\hline
\end{tabular}
\end{subtable}
\caption{Expected cost of abandonments for the parameter values given in Table $\ref{langt2}$ for English, Hindi and Mandarin speaking customers separately.}
    \label{langeg2t}
\end{table}

We observe that the English speaking customers abandon the least when there is no reservation in place, while Hindi and Mandarin speaking customers abandon the least when their respective agents are fully reserved. The expected number of abandonments for Hindi speaking customers increases as the number of reserved Mandarin speaking agents increases, since more English queries are then routed to Hindi speaking agents -- and vice versa. if all Hindi speaking agents are reserved, the number of abandonments for Hindi speaking customers becomes independent of the Mandarin agent reservation level, and the same holds in the reverse scenario.

Additionally, the expected number of abandonments for Hindi speaking customers is higher than that for Mandarin speaking customers. This is due to the asymmetric call assignment rule, under which, when the number of available Hindi and Mandarin speaking agents above their thresholds is equal, the English call is assigned to a Hindi speaking agent.

\section{Discussion}\label{dis}

The methods presented in this paper have so far been discussed in the context of a single time interval. However, in practical applications, these methods would be implemented sequentially over multiple intervals. Although the examples above begin with empty systems, the analysis is equally valid for any initial state. In a real-world call centre, one could begin with an empty system, determine the optimal policy for a given time interval, and then reapply the method at the end of that interval using the system’s current state as the new initial state. This approach captures the time dependency of the optimal policy. Additionally, when applying Bellman’s equation, one could experiment with different boundary conditions and select the one most appropriate for the specific operational objective, for example, a boundary value function proportional to the number of customers waiting in the queue at time step $M$.

Many call centre managers think that it is important to serve customers according to a first-come-first-served (FCFS) discipline. We do not assume any specific discipline in the analysis presented in this paper. A call centre implementing the proposed methods may choose to follow FCFS discipline within each customer level. However, the class of reservation policies is not FCFS across different levels.

Even with the simplifications that occur when we consider reservation policies, as opposed to general policies, our analysis is limited by the size of the system. The larger the size of the system, the bigger is the matrix inversion step to calculate the Laplace transforms, and hence the computational expense increases. But, unlike using Bellman's equation, changing $t$ does not increase the computational expense. The time to find the optimal policy will also depend on the number of possible values of $\mathbf{\Theta}$. 

For $\ell = 10$, to calculate the expected cost of abandonments for a given $\mathbf{\Theta}$, it took $300$ seconds for the hierarchical model with four levels -- examples \ref{eg1}, \ref{eg2} and \ref{eg3}. Using $\ell > 15$ for the hierarchical model with four levels led to memory allocation issues on our laptop. In the multilingual call centre model -- example \ref{eg4}, for $\ell = 12$, it took $2-3$ seconds and for $\ell = 20$, it took $400$ seconds to calculate the expected cost of abandonments for a given $\Theta$. This is because the multilingual call centre model has a lower dimensional state space. 

We also employed Black-box optimisation in Julia with \textit{MaxSteps} parameters equal to $3$ and it took $13$ seconds to find the best policy for $\ell = 10$ for the multilingual call centre model. The result is not guaranteed to be optimal, but comparing the Black-box optimisation results with the enumeration results showed that it was able to find the optimal policy in the first iteration. For $\ell = 20$, it took $2937$ seconds to perform three steps of the Black box iteration. However, it had already found the optimal value in $927$ seconds in the first step.

\section{Conclusion}\label{con}

In this paper, we have used different methods to find the optimal finite-horizon call allocation policy for small call centres. Backward induction and an infinite time model with discounting are effective methods for finding the class of policies that are likely to be optimal, but they do not calculate the exact expected cost of abandonments in continuous time. 

The method introduced above employs first-step analysis, Laplace transforms and numerical inversion. In this paper we have considered the expected cost of abandonments. In a forthcoming paper, we will describe how an extension of of these methods can be used to find the waiting time distributions. 

We have an ongoing objective to increase the computational efficiency of the method so that it can be applied to larger systems.

\textbf{Acknowledgements.} This research was funded by the Australian Government through the Australian Research Council Industrial Transformation Training Centre in Optimisation Technologies, Integrated Methodologies, and Applications (OPTIMA), Project ID IC200100009. It was also partially funded by Probe CX and supported by The University of Melbourne’s Research Computing Services and the Petascale Campus Initiative. The second author is also funded by the Melbourne Research Scholarship.

\textbf{Competing interests declaration.} The authors declare none.
\newpage

\printbibliography[title= {References}]

\appendix

\section{Transition rates for the hierarchical model with $\nu = 4$}\label{SubS:Transition}

Transitions between states occur when there is an arrival, a service or an abandonment at one of the levels 1 to 4. Below is a list of the transitions that ensue when each of these events occurs when the system is in state $(a,a_1,b,b_1,c,c_1,d)$. 

\begin{itemize}
    \item If a level 1 caller arrives, there are three possibilities.
    \begin{itemize}
        \item The caller is allocated to a level 1 agent if at least one of them is available leading to the state $(a+1,a_1,b,b_1,c,c_1,d)$.
        \item The caller is allocated to a level 2 agent (if no level 1 agent is available, no level 2 call is in waiting, and there are more than $\Theta_2$ level 2 agents available) leading to the state $(a+1,a_1+1,b,b_1,c,c_1,d)$.
        \item The caller joins the queue of level one customers if neither of the first two conditions are satisfied leading to the state $(a+1,a_1,b,b_1,c,c_1,d)$.
    \end{itemize} 
    Hence, there is a transition to the state $(a+1,a_1+\mathbbm{1}(a-a_1 \geq k_1, b - b_1 < k_2 - \Theta_2 - a_1),b,b_1,c,c_1,d)$ with rate $\lambda_1$ if $a+b+c+d < \ell$, otherwise the caller is lost, costing $\beta$.
    \item Similarly, a transition to the state $(a,a_1,b+1,b_1+\mathbbm{1}(b-b_1 \geq k_2, c - c_1 < k_3 - \Theta_3 - b_1),c,c_1,d)$ happens with rate $\lambda_2$ and a transition to the state $(a,a_1,b,b_1,c+1,c_1+\mathbbm{1}(c-c_1 \geq k_3, d < k_4 - \Theta_4 - c_1),d)$ happens with the rate $\lambda_3$ if $a+b+c+d < \ell$.
    \item When a level 4 caller arrives, there is only one possibility, that is, a transition to the state $(a,a_1,b,b_1,c,c_1,d+1)$ with rate $\lambda_4$ if $a+b+c+d < \ell$.
    \item A service of level 1 caller by a level 1 agent happens with rate $\min(k_1,a-a_1)\mu_1$ leading to the state $(a-1,a_1,b,b_1,c,c_1,d)$.
    \item When there is a service of a level 1 caller by a level 2 agent, there are two possibilities. 
    \begin{itemize}
        \item That agent becomes free or starts serving a level 2 caller if there are any in waiting leading to the state $(a-1,a_1-1,b,b_1,c,c_1,d)$.
        \item The agent starts serving another level 1 caller leading to the state $(a-1,a_1,b,b_1,c,c_1,d)$.
    \end{itemize}
    Hence, there is transition to the state $(a-1,a_1 - (1-\mathbbm{1}(a-a_1 > k_1, b - b_1 = k_2 - \Theta_2 - a_1),b,b_1,c,c_1,d)$ with rate $a_1 \mu_1'$.
    \item Similarly, when there is a service of a level 2 caller by a level 2 agent, that agent can either serve another level 2 caller, serve a level 1 caller or stay available. The system transitions to the state $(a,a_1+\mathbbm{1}(a-a_1 > k_1, b - b_1 = k_2 - \Theta_2 - a_1),b-1,b_1,c,c_1,d)$ with rate $\min(b-b_1,k_2-a_1)\mu_2$.
    \item Similarly, the system transitions to 
    \begin{itemize}
        \item $(a,a_1,b-1,b_1 - (1-\mathbbm{1}(b-b_1 > k_2 - a_1, c = k_3 - \Theta_3 - b_1),c,c_1,d)$ with rate $b_1 \mu_2'$ when there is a service of level 2 caller by a level 3 agent,
        \item $(a,a_1,b,b_1+\mathbbm{1}(b-b_1 > k_2, c - c_1 = k_3 - \Theta_3 - b_1),c-1,c_1,d)$ with rate $\min(c-c_1,k_3-b_1)\mu_3$ when there is a service of a level 3 caller by a level 3 agent,
        \item $(a,a_1,b,b_1,c-1,c_1 - (1-\mathbbm{1}(c-c_1 > k_3 - b_1, d = k_4 - \Theta_4 - c_1),d)$ with rate $c_1 \mu_3'$ when a level 3 caller is served by a level 4 agent and
        \item $(a,a_1,b,b_1,c,c_1+\mathbbm{1}(c-c_1 > k_3 - b_1, d = k_4 - \Theta_4 - c_1),d-1)$ with rate $\min(d,k_4-c_1)\mu_4$ when there is a service of level 4 caller.
    \end{itemize}
    \item The transitions in case of abandonments are straightforward. There is a transition to 
    \begin{itemize}
        \item $(a-1,a_1,b,b_1,c,c_1,d)$ when a level 1 caller abandons with rate $\max(0,a-(k_1 + a_1))\theta_1$ costing $\gamma_1$,
        \item $(a,a_1,b-1,b_1,c,c_1,d)$ when a level 2 caller abandons with rate $\max(0,b-(k_2-a_1+b_1))\theta_2$ costing $\gamma_2$,
        \item $(a,a_1,b,b_1,c-1,c_1,d)$ when a level 3 caller abandons with rate $\max(0,c-(k_3-b_1+c_1))\theta_3$ costing $\gamma_3$ and 
        \item $(a,a_1,b,b_1,c,c_1,d-1)$ when a level 4 caller abandons with rate $\max(0,d-(k_4-c_1))\theta_4$ costing $\gamma_4$.
    \end{itemize}
\end{itemize}

\newpage
\section{Laplace transform expressions for the expected cost of losses and abandonments for the hierarchical model with $\nu = 4$}\label{laplace}

The Laplace transform
\begin{equation}
\Tilde{C}_{a,a_1,b,b_1,c,c_1,d}(\mathbf{\Theta},s) = \dfrac{\xi_{a,a_1,b,b_1,c,c_1,d}(\mathbf{\Theta},s)}{\psi_{a,a_1,b,b_1,c,c_1,d}(s)}
\end{equation}
where
\begin{equation}
    \begin{split}
        \xi_{a,a_1,b,b_1,c,c_1,d}(\mathbf{\Theta},s) = &(\lambda_1 \Tilde{C}_{a+1,a_1+\mathbbm{1}(a-a_1 \geq k_1, b - b_1 < k_2 - \Theta_2 - a_1),b,b_1,c,c_1,d}(\mathbf{\Theta},s) \\&+ \lambda_2 \Tilde{C}_{a,a_1,b+1,b_1+\mathbbm{1}(b-b_1 \geq k_2, c < k_3 - \Theta_3 - b_1),c,c_1,d}(\mathbf{\Theta},s)\\&+ \lambda_3 \Tilde{C}_{a,a_1,b,b_1,c+1,c_1+\mathbbm{1}(c-c_1 \geq k_3, d < k_4 - \Theta_4 - c_1),d}(\mathbf{\Theta},s)\\ & + \lambda_4 \Tilde{C}_{a,a_1,b,b_1,c,c_1,d+1}(\mathbf{\Theta},s)) \mathbbm{1}_{a+b+c+d < \ell}+\min(k_1,a-a_1)\mu_1 \Tilde{C}_{a-1,a_1,b,b_1,c,c_1,d}(\mathbf{\Theta},s) \\ &+ a_1 \mu_1' \Tilde{C}_{a-1,a_1 - (1-\mathbbm{1}(a-a_1 > k_1, b - b_1 = k_2 - \Theta_2 - a_1),b,b_1,c,c_1,d}(\mathbf{\Theta},s))\\ &+\min(b-b_1,k_2-a_1)\mu_2 \Tilde{C}_{a,a_1+\mathbbm{1}(a-a_1 > k_1, b - b_1 = k_2 - \Theta_2 - a_1),b-1,b_1,c,c_1,d}(\mathbf{\Theta},s)\\ & + b_1 \mu_2' \Tilde{C}_{a,a_1,b-1,b_1 - (1-\mathbbm{1}(b-b_1 > k_2 - a_1, c = k_3 - \Theta_3 - b_1),c,c_1,d}(\mathbf{\Theta},s)) \\&+\min(c-c_1,k_3-b_1)\mu_3 \Tilde{C}_{a,a_1,b,b_1+\mathbbm{1}(b-b_1 > k_2, c - c_1 = k_3 - \Theta_3 - b_1),c-1,c_1,d}(\mathbf{\Theta},s)\\ & + c_1 \mu_3' \Tilde{C}_{a,a_1,b,b_1,c-1,c_1 - (1-\mathbbm{1}(c-c_1 > k_3 - b_1, d = k_4 - \Theta_4 - c_1),d}(\mathbf{\Theta},s)) \\ &+ \min(d,k_4-c_1)\mu_4 \Tilde{C}_{a,a_1,b,b_1,c,c_1+\mathbbm{1}(c-c_1 > k_3 - b_1, d = k_4 - \Theta_4 - c_1),d-1}(\mathbf{\Theta},s) \\ &+\max(0,a-(k_1 + a_1))\theta_1(\Tilde{C}_{a-1,a_1,b,b_1,c,c_1,d}(\mathbf{\Theta},s) + \gamma_1/s)\\ &+\max(0,b-(k_2-a_1+b_1))\theta_2(\Tilde{C}_{a,a_1,b-1,b_1,c,c_1,d}(\mathbf{\Theta},s)+\gamma_2/s)\\ &+\max(0,c-(k_3-b_1+c_1))\theta_3(\Tilde{C}_{a,a_1,b,b_1,c-1,c_1,d}(\mathbf{\Theta},s)+\gamma_3/s)\\ &+\max(0,d-(k_4-c_1))\theta_4(\Tilde{C}_{a,a_1,b,b_1,c,c_1,d-1}(\mathbf{\Theta},s)+\gamma_4/s)\\ &+(\lambda_1+\lambda_2+\lambda_3+\lambda_4)\beta\mathbbm{1}_{a+b+c+d = \ell}/s,
    \end{split}
\end{equation}
\begin{equation}
    \begin{split}
        \psi_{a,a_1,b,b_1,c,c_1,d}(s) = &(\lambda_1+\lambda_2+\lambda_3+\lambda_4)\mathbbm{1}_{a+b+c+d < \ell} + \min(k_1,a-a_1)\mu_1 + a_1 \mu_1' + \min(b-b_1,k_2-a_1)\mu_2 \\ & + b_1 \mu_2' + \min(c-c_1,k_3-b_1)\mu_3 + c_1 \mu_3' + \min(d,k_4-c_1)\mu_4 + \max(0,a-(k_1 + a_1))\theta_1 \\ &+ \max(0,b-(k_2-a_1+b_1))\theta_2 + \max(0,c-(k_3-b_1+c_1))\theta_3\\ &+\max(0,d-(k_4-c_1))\theta_4+s
    \end{split}
\end{equation}
and
$\Tilde{C}_{a,a_1,b,b_1,c,c_1,d}(\mathbf{\Theta},s) = 0$ for $(a,a_1,b,b_1,c,c_1,d) \notin S$. 

\end{document}